\documentclass[11pt]{article}
\usepackage{graphicx}
\usepackage{amsmath,amsfonts,amssymb}
\usepackage{citesort}
\usepackage{url}
\usepackage{amsmath}
\usepackage{graphicx}
\usepackage{epsfig}
\usepackage{epstopdf}
\usepackage{color}
\def\tto{\;{\lower 1pt \hbox{$\rightarrow$}}\kern -10pt
\hbox{\raise 2pt \hbox{$\rightarrow$}}\;}
\def\Hat{\widehat}

\def\Bar{\overline}
\def\ra{\rangle}
\def\la{\langle}

\def\B{I\!\!B}
\def\h{\hfill\Box}
\def\R{\Bbb R}
\def\N{\Bbb N}
\def\ox{\bar{x}}

\def\h{\hfill\square}

\def\O{\Omega}
\def\ph{\varphi}

\def\oR{\Bar{\R}}

\newcounter{lk}

\begin{document}

\begin{center}
\vspace*{0.3in} {\bf SUBDIFFERENTIAL FORMULAS FOR A CLASS OF NONCONVEX INFIMAL CONVOLUTIONS}\\[2ex]
Nguyen Mau Nam\footnote{Fariborz Maseeh Department of
Mathematics and Statistics, Portland State University, PO Box 751, Portland, OR 97207, United States (mau.nam.nguyen@pdx.edu). The research of Nguyen Mau Nam was partially supported by
the Simons Foundation under grant \#208785.}
\end{center}
\small{\bf Abstract:} In this paper, we provide a number of subdifferential formulas for a class of nonconvex infimal convolutions in normed spaces. The formulas obtained unify several results on subdifferentials of the distance function and the minimal time function. In particular, we generalize and validate the results obtained recently by Zhang, He, and Jiang \cite{ZHJ}.\\[1ex]
\vspace*{0,05in} \noindent {\bf Key words.}  subdifferentials; infimal convolution; minimal time function

\newtheorem{Theorem}{Theorem}[section]
\newtheorem{Proposition}[Theorem]{Proposition}
\newtheorem{Remark}[Theorem]{Remark}
\newtheorem{Lemma}[Theorem]{Lemma}
\newtheorem{Corollary}[Theorem]{Corollary}
\newtheorem{Definition}[Theorem]{Definition}
\newtheorem{Example}[Theorem]{Example}
\renewcommand{\theequation}{\thesection.\arabic{equation}}
\normalsize

\section{Introduction}

Let $X$ be a normed space and let $\O$ be a nonempty subset of $X$. The \emph{distance function} to $\O$ is defined on $X$ by
\begin{equation}\label{ds}
d(x;\O):=\inf\big\{\|x-\omega\|\; \big|\; \omega\in \O\big\}.
\end{equation}
The distance function (\ref{ds}) belongs to a larger class of functions called the \emph{minimal time function} defined in what follows. Given a nonempty closed bounded convex set $F$, define the \emph{Minkowski gauge} associated with $F$ by
\begin{equation}\label{mg}
\rho_F(x):=\inf\big\{t\geq 0\; \big|\; x\in tF\big\}.
\end{equation}
The Minkwoski gauge (\ref{mg}) reduces to the normed function when $F$ is the closed unit ball of $X$. Based on the Minkowski gauge, the minimal time function to the set $\O$ is defined by
\begin{equation}\label{mt}
\mathcal{T}_F(x;\O):=\inf\big\{\rho_F(\omega-x)\; \big|\; \omega\in \O\big\},
\end{equation}
which is obviously a more general form of the distance function (\ref{ds}).

The minimal time function forms an interesting class of nonsmooth functions due to its intrinsic nondifferentiability. Subdifferential formulas for this class functions in both convex and nonconvex settings have been of great interest in the literature; see \cite{cowo,JH,HN,MN10,MN,NZ} and the references therein. It is well known that the subdifferential in the sense of convex analysis of the distance function (\ref{ds}) can be computed using the following \emph{infimal convolution} representation:
\begin{equation}\label{dr}
d(x;\O)=\inf\{\|y-x\|+\delta_\O(y)\; |\; y\in X\},
\end{equation}
where $\delta_\O$ is the \emph{indicator function} associated with $\O$ given by $\delta(x;\O)=0$ if $x\in \O$, and $\delta(x;\O)=\infty$ otherwise. However, a similar approach for nonconvex setting has not been available due to the lack of subdifferential formulas for nonconvex infimal convolutions.

Let $X^*$ denote the topological dual of $X$ and let $\la \cdot,\cdot\ra$ denote the paring between $X$ and $X^*$. Given a convex function $g: X\to \oR$ and given $\ox\in \mbox{\rm dom}\,g:=\{x\in X\; |\; |g(x)|<\infty\}$, the subdifferential of $g$ in the sense of convex analysis at $\ox$ is defined by
\begin{equation*}
\partial g(\ox):=\big\{x^*\in X^*\; \big |\; \la x^*, x-\ox\ra\leq g(x)-g(\ox)\; \mbox{\rm for all }x\in X\big\}.
\end{equation*}
In the same setting, but the convexity of the function $g$ is not assumed, the \emph{$\epsilon-$Fr\'echet subdifferential} ($\epsilon\geq 0$) of $g$ at $\ox$ is the set
\begin{equation*}
\Hat\partial_\epsilon g(\ox):=\big\{x^*\in X^*\; \big |\; \liminf_{x\to \ox}\dfrac{g(x)-g(\ox)-\la x^*, x-\ox\ra}{\|x-\ox\|}\geq -\epsilon\big\}.
\end{equation*}
In the case where $\epsilon=0$, the set $\Hat\partial_0g(\ox)$ is called the \emph{Fr\'echet subdifferential} of $g$ at $\ox$ and is denoted simply by $\Hat\partial g(\ox)$.

Another useful subdifferential construction called \emph{proximal subdifferential} is defined by
\begin{equation*}
\partial_Pg(\ox):=\big\{x^*\in X^*\; \big |\; \liminf_{x\to \ox}\dfrac{g(x)-g(\ox)-\la x^*, x-\ox\ra}{\|x-\ox\|^2}> -\infty\big\}.
\end{equation*}
If $g$ is convex, both Fr\'echet and proximal subdifferential constructions reduce to the subdifferential in the sense of convex analysis, while $\Hat\partial_\epsilon g(\ox)$ has the following representation:
\begin{equation*}
\Hat\partial_\epsilon g(\ox)=\big\{x^*\in X^*\; \big |\; \la x^*, x-\ox\ra \leq g(x)-g(\ox)+\epsilon \|x-\ox\|\; \mbox{\rm for all }x\in X\big\}.
\end{equation*}
The readers are referred to the monographs \cite{C,CLSW,mor} for more properties of these subdifferential constructions as well as their applications.

In a recent paper published in \emph{Optimization Letters}, Zhang, He, and Jiang \cite{ZHJ} introduced and studied the so-called \emph{perturbed minimal time function} $T^f_F: X\to \oR:=[-\infty, \infty]$ defined by
\begin{equation}\label{pm}
T^f_F(x):=\inf\big\{T(x,y)+f(y)\; \big |\; y\in X\big\},
\end{equation}
where $T(x,y)=\inf\{t\geq 0\; |\; y-x\in t F\}=\rho_F(y-x)$, and $f: X\to (-\infty, \infty]$ is an extended-real-valued function. Then subdifferential formulas of Fr\'echet and proximal types were developed for this class of functions at points belonging to the set
\begin{equation}\label{ss}
S_0:=\big\{x\in \mbox{\rm dom}\,(T^f)\; \big |\; T^f_F(x)=f(x)\big\}.
\end{equation}
Note that in the setting of (\ref{pm}) with $f(x)=\delta_\O(x)$ and $T(x,y)=\|y-x\|$, it is obvious that $S_0=\O$.

In this paper, we consider a general class of nonconvex infimal convolutions and investigate its subdifferential properties.  Let $\ph: X\to (-\infty,\infty]$ and let $f: X\to (-\infty,\infty]$ be extended-real-valued functions. Consider the \emph{infimal convolution} of $\ph$ and $f$:
\begin{equation}\label{ic}
T_\ph^f(x):=\inf\{\ph(y-x)+f(y)\; |\; y\in X\}.
\end{equation}
Define
\begin{equation}\label{s}
S_0:=\{x\in \mbox{\rm dom}\,(T^f_\ph)\; |\; T_\ph^f(x)=f(x)\}.
\end{equation}
Note that the infimal convolution (\ref{ic}) also covers another class of functions called the \emph{perturbed distance function} $d_\O^J: X\to \oR$ given by
\begin{equation}\label{pd}
d_\O^J(x):=\inf\big\{\|y-x\|+J(y)\; \big|\; y\in \O\big\}
\end{equation}
where $J: X\to (-\infty, \infty]$ is an extended-real-valued function. In this setting of the infimal convolution (\ref{ic}) with $\ph(x)=\|x\|$, $f(x)=J(x)+\delta_\O(x)$, we obtain the perturbed distance function (\ref{pd}) for which the set $S_0$ reduces to
\begin{equation}\label{pds}
S_0=\big\{x\in \mbox{\rm dom}\,(d_\O^J)\; \big|\; d_\O^J(x)=J(x)\big\}.
\end{equation}
Subdifferential properties of the perturbed distance function were the topics of study in \cite{MLY,WLX}.

Our paper is organized as follows. In Section 2, we develop $\epsilon-$Fr\'echet subdifferential formulas for the infimal convolution (\ref{s}) at points belonging to the set $S_0$. Section 3 is devoted to corresponding  \emph{H\"{o}lder subdifferential} formulas. The results we obtain unify many related results for the distance function and the minimal time function available in the literature. In particular, we generalize and validate the results obtained recently by Zhang, He, and Jiang in \cite{ZHJ}.

\section{Fr\'echet Subdifferential Formulas}

This section focuses on $\epsilon-$Fr\'echet subdifferential formulas for the infimal convolution (\ref{ic}) and its specifications. In the proposition below, we give an upper estimate for the $\epsilon-$Fr\'echet subdifferential of the infimal convolution (\ref{ic}). Note that this result is well known in the convex case, but it is new in the nonconvex case.
\begin{Proposition}\label{p1} Consider the infimal convolution {\rm (\ref{ic})} and the set $S_0$ given by {\rm (\ref{s})} with $\ox\in S_0$. Suppose that $\ph(0)=0$. Given $\epsilon\geq 0$, one has
\begin{equation*}
\Hat\partial_\epsilon T_\ph^f(\ox)\subset \Hat\partial_\epsilon f(\ox)\cap \big[-\Hat\partial_\epsilon \ph(0)\big].
\end{equation*}
\end{Proposition}
{\bf Proof.} Fix any $x^*\in \Hat\partial_\epsilon T^f_\ph(\ox)$. Then for any $\eta>0$ there exists $\delta>0$ such that
\begin{equation*}
\la x^*, x-\ox\ra \leq T_\ph^f(x)-T_\ph^f(\ox)+(\epsilon+\eta) \|x-\ox\|\; \mbox{\rm whenever }\|x-\ox\|<\delta.
\end{equation*}
Since $T_\ph^f(\ox)=f(\ox)$, it follows that
\begin{equation*}
\la x^*, x-\ox\ra \leq T_\ph^f(x)-f(\ox)+(\epsilon+\eta) \|x-\ox\|\; \mbox{\rm whenever }\|x-\ox\|<\delta.
\end{equation*}
From the definition of $T_\ph^f$, it is obvious that $T_\ph^f(x)\leq \ph(x-x)+f(x)=\ph(0)+f(x)=f(x)$ for all $x\in X$. Thus,
\begin{equation*}
\la x^*, x-\ox\ra \leq f(x)-f(\ox)+(\epsilon+\eta) \|x-\ox\|\; \mbox{\rm whenever }\|x-\ox\|<\delta,
\end{equation*}
which implies that $x^*\in \Hat\partial_\epsilon f(\ox)$.

Fix any $v\in X$ with $\|v\|<\delta$. Then $\ox-v\in \B(\ox; \delta)$ and
\begin{equation*}
T^f_\ph(\ox-v)=\inf\big\{\ph(y-(\ox-v))+f(y)\; \big |\; y\in X\big\}\leq \ph(\ox-(\ox-v))+f(\ox)=\ph(v)+f(\ox),
\end{equation*}
which implies
\begin{equation*}
\la x^*, -v\ra \leq T_\ph^f(\ox-v)-f(\ox)+(\epsilon+\eta) \|v\|\leq \ph(v)+(\epsilon+\eta) \|v\|.
\end{equation*}
Since $\eta>0$ is arbitrary, it follows that $-x^*\in \Hat\partial_\epsilon \ph(0)$, and hence $x^*\in -\Hat\partial_\epsilon \ph(0)$. $\h$

Let $g: X\to (-\infty, \infty]$ be an extended-real-valued function. We say that $g$ is \emph{coercive} with constant $m>0$ on $X$ if
\begin{equation*}
m\|x\|\leq g(x)\; \mbox{\rm for all }x\in X.
\end{equation*}
We also say that $g$ satisfies a \emph{center-Lipschitz/calm} condition on a set $D\subset X$ at $\ox\in D$ with constant $\ell\geq 0$ if
\begin{equation*}
|g(x)-g(\ox)|\leq \ell \|x-\ox\|\; \mbox{\rm for all }x\in D.
\end{equation*}
In the proposition below, we prove that the Minkowski gauge (\ref{mg}) is coercive. For the convenience of representation, we assume that $F$ is nonzero. In the case where $F=\{0\}$, we can easily verify that the Minkowski gauge (\ref{mg}) is also coercive with constant $m$, where $m$ is any positive real number.
\begin{Proposition}\label{p2} Suppose that $F$ is a nonempty closed bounded convex set that is nonzero. Then the Minkowski gauge {\rm (\ref{mg})} is subadditive, positively homogeneous with $\rho_F(0)=0$, and coercive with constant $m:=\|F\|^{-1}$, where
\begin{equation*}
\|F\|:=\sup\big\{\|u\|\; \big |\; u\in F\big\}.
\end{equation*}
Moreover,
\begin{equation}\label{mgs}
\partial\rho_F(0)=\big\{x^*\in X^*\; \big |\; \sup_{u\in F}\la x^*,u\ra \leq 1\big\}.
\end{equation}
\end{Proposition}
{\bf Proof.} It follows from the definition that $\rho_F$ is subadditive and positively homogeneous with $\rho_F(0)=0$, so it is a convex function. Obviously, $\rho_F(x)\geq m\|x\|$ if $x\notin \mbox{\rm dom}\,\rho_F$. Fix any $x\in \mbox{\rm dom}\,\rho_F$ and let $(t_k)$ be a sequence of nonnegative numbers such that $t_k\to \rho_F(x)$ and $x\in t_k F$ for every $k$. Then
\begin{equation*}
\|x\|\leq t_k\|F\|.
\end{equation*}
It follows by passing to a limit as $k\to \infty$ that $\|F\|^{-1}\|x\|\leq \rho_F(x)$, which justifies the coercivity of $\rho_F$.

The subdifferential formula (\ref{mgs}) is well known, but we provide a proof for the convenience of the readers. From the definition, one has that $\rho_F(0)=0$. Assuming that $x^*\in \partial \rho_F(0)$ implies
\begin{equation*}
\la x^*, u\ra \leq \rho_F(u)\leq 1\; \mbox{\rm for all }u\in F,
\end{equation*}
which verifies the inclusion $\subset$ in (\ref{mgs}). Now suppose that $\sup_{u\in F}\la x^*,u\ra \leq 1$ and fix any $x\in \mbox{\rm dom}\,\rho_F$. Let $(t_k)$ be a sequence of nonnegative numbers such that $t_k\to \rho_F(x)$ and $x\in t_k F$ for every $k$. For every $k\in \N$, find $u_k\in F$ with $x=t_ku_k$.  Then
\begin{equation*}
\la x^*, x\ra =\lim_{k\to \infty}\la x^*, t_ku_k\ra =\lim_{k\to \infty}t_k\la x^*, u_k\ra\leq \lim_{k\to \infty}t_k=\rho_F(x).
\end{equation*}
This implies $x^*\in \partial\rho_F(0)$ and completes the proof. $\h$

The theorem below generalizes the result of \cite[Theorem~3.1]{ZHJ} from the perturbed minimal time function (\ref{pm}) to the general infimal convolution (\ref{ic}). Note that we only assume that $f$ satisfies a center-Lipschitz condition on its domain instead of the whole space $X$ as in \cite{ZHJ}. This is important because the indicator function $\delta_\O$ obviously satisfies a center-Lipschitz condition on its domain $\O$ with constant $\ell=0$, but it does not satisfies a center-Lipschitz condition on $X$.

\begin{Theorem}\label{t1} Consider the infimal convolution {\rm (\ref{ic})} in which $\ph(0)=0$ and consider the set $S_0$ given by {\rm (\ref{s})} with $\ox\in S_0$. Suppose that and $\ph$ is coercive on $X$ with constant $m>0$ and $f$ satisfies a center-Lipschitz condition on $D:=\mbox{\rm dom}\,f$ with constant $\ell$ where $0\leq \ell <m$. Given $\epsilon\geq 0$ and $x^*\in \Hat\partial_\epsilon f(\ox)\cap \big[-\Hat\partial_\epsilon \ph(0)\big]$, one has
\begin{equation}\label{ff1}
x^*\in \Hat\partial_{\alpha\epsilon}T_\ph^f(\ox), \mbox{\rm where }\alpha:=2(\|x^*\|+m)(m-\ell)^{-1}+1.
\end{equation}
Moreover,
\begin{equation}\label{f1}
\Hat\partial T_\ph^f(\ox)= \Hat\partial f(\ox)\cap \big[-\Hat\partial \ph(0)\big].
\end{equation}
\end{Theorem}
{\bf Proof.} Suppose by contradiction that there exists $x^*\in \Hat\partial_\epsilon f(\ox)\cap \big[-\Hat\partial_\epsilon \ph(0)\big]$, but (\ref{ff1}) is not satisfied. Then there exist $\sigma>0$ and a sequence $(x_k)$ that converges to $\ox$ such that
\begin{equation}\label{e1}
T_\ph^f(x_k)<T_\ph^f(\ox)+\la x^*, x_k-\ox\ra -(\alpha\epsilon+\sigma) \|x_k-\ox\|\; \mbox{\rm for every }k,
\end{equation}
which implies
\begin{equation*}
T_\ph^f(x_k)<T_\ph^f(\ox)+\la x^*, x_k-\ox\ra\leq f(\ox)+\|x^*\|\|x_k-\ox\|.
\end{equation*}
Observe also that $x_k\neq \ox$ for every $k$. From the definition of $T_\ph^f(x_k)$, find $y_k\in X$ such that
\begin{equation*}
\ph(y_k-x_k)+f(y_k)\leq T_\ph^f(x_k)+\|x_k-\ox\|^2.
\end{equation*}
It follows that
\begin{equation*}
\ph(y_k-x_k)+f(y_k)\leq T_\ph^f(x_k)+\|x_k-\ox\|^2\leq f(\ox)+\|x^*\|\|x_k-\ox\|+\|x_k-\ox\|^2.
\end{equation*}
Then $y_k\in \mbox{\rm dom}\,f$ for every $k$ and
\begin{align*}
m\|y_k-\ox\|&\leq m\|y_k-x_k\|+m \|x_k-\ox\|\\
&\leq \ph(y_k-x_k)+m \|x_k-\ox\|\leq f(\ox)-f(y_k)+(\|x^*\|+m)\|x_k-\ox\|+\|x_k-\ox\|^2\\
&\leq \ell \|y_k-\ox\|+(\|x^*\|+m)\|x_k-\ox\|+\|x_k-\ox\|^2.
\end{align*}
It follows that
\begin{equation*}
(m-\ell)\|y_k-\ox\|\leq (\|x^*\|+m)\|x_k-\ox\|+\|x_k-\ox\|^2,
\end{equation*}
which implies
\begin{equation*}
\|y_k-\ox\|\leq (m-\ell)^{-1}(\|x^*\|+m)\|x_k-\ox\|+o(\|x_k-\ox\|)\to 0\; \mbox{\rm as }k\to \infty.
\end{equation*}
Since $x^*\in \Hat\partial_\epsilon f(\ox)$, given any $\eta>0$, find $\delta>0$ such that
\begin{equation}\label{fr1}
\la x^*, x-\ox\ra \leq f(x)-f(\ox)+(\epsilon+\eta) \|x-\ox\|\; \mbox{\rm whenever }\|x-\ox\|<\delta.
\end{equation}
We can assume without loss of generality that $\|y_k-\ox\|<\delta$ for every $k$, and hence (\ref{fr1}) holds with $x:=y_k$. Taking into account that $x^*\in -\Hat\partial_\epsilon \ph(0)$ and that $\|y_k-x_k\|\to 0$ as $k\to \infty$, we can also assume without loss of generality that
\begin{equation*}
-(\epsilon+\eta) \|y_k-x_k\|\leq \ph(y_k-x_k)-\la x^*, x_k-y_k\ra.
\end{equation*}
Then the following estimates hold:
\begin{align*}
T_\ph^f(x_k)-T_\ph^f(\ox)-\la x^*, x_k-\ox\ra &=T_\ph^f(x_k)-f(\ox)-\la x^*, x_k-y_k\ra -\la x^*, y_k-\ox\ra\\
&\geq T_\ph^f(x_k)-f(\ox)-\la x^*, x_k-y_k\ra-\big[f(y_k)-f(\ox)+(\epsilon+\eta) \|y_k-\ox\|\big]\\
&= T_\ph^f(x_k)-f(y_k)-\la x^*, x_k-y_k\ra-(\epsilon+\eta) \|y_k-\ox\|\\
&\geq \ph (y_k-x_k)-\|x_k-\ox\|^2-\la x^*, x_k-y_k\ra-(\epsilon+\eta) \|y_k-\ox\|\\
&\geq -(\epsilon+\eta)\|y_k-x_k\|-\|x_k-\ox\|^2-(\epsilon+\eta) \|y_k-\ox\|\\
&\geq -(\epsilon+\eta)\|y_k-\ox\|-(\epsilon+\eta)\|x_k-\ox\|-\|x_k-\ox\|^2-(\epsilon+\eta) \|y_k-\ox\|\\
&\geq -(2\epsilon+2\eta)\|y_k-\ox\|-(\epsilon+\eta)\|x_k-\ox\|-\|x_k-\ox\|^2\\
&\geq -(2\epsilon+2\eta) (m-\ell)^{-1}(\|x^*\|+m)\|x_k-\ox\|-(\epsilon+\eta)\|x_k-\ox\|-o(\|x_k-\ox\|).
\end{align*}
Comparing with (\ref{e1}) yields
\begin{equation*}
-(2\epsilon+2\eta) (m-\ell)^{-1}(\|x^*\|+m)\|x_k-\ox\|-(\epsilon+\eta)\|x_k-\ox\|-o(\|x_k-\ox\|)\leq -(\alpha\epsilon+\sigma) \|x_k-\ox\|,
\end{equation*}
which implies
\begin{equation*}
\alpha\epsilon+\sigma\leq (2\epsilon+2\eta) (m-\ell)^{-1}(\|x^*\|+m)+\epsilon+\eta.
\end{equation*}
Letting $\eta\to 0^+$, one has that $$\alpha\epsilon+\sigma\leq 2\epsilon(m-\ell)^{-1}(\|x^*\|+m)+\epsilon=\alpha\epsilon,$$
which is a contradiction. We have proved the first statement.

The subdifferential equality (\ref{f1}) follows from the first statement and Proposition \ref{p1} with $\epsilon=0$. The proof is now complete. $\h$

As a corollary, we obtain \cite[Theorem~3.1]{ZHJ} with some validation.
\begin{Corollary}\label{c1}
Consider the infimal convolution {\rm (\ref{pm})} and the set $S_0$ given by {\rm (\ref{ss})} with $\ox\in S_0$. Suppose that $f$ satisfies a center-Lipschitz condition on $D:=\mbox{\rm dom}\,f$ with constant $\ell$ where $0\leq \ell <\|F\|^{-1}$. Then
\begin{equation}\label{f2}
\Hat\partial T^f_F(\ox)= \Hat\partial f(\ox)\cap \big\{x^*\in X^*\; \big |\; \sup_{u\in F}\la -x^*, u\ra\leq 1\big\}.
\end{equation}
\end{Corollary}
{\bf Proof.} The subdifferential formula (\ref{f2}) follows from Theorem \ref{t1} and Proposition \ref{p2}. $\h$

Using Corollary \ref{c1} and the fact that the norm function is coercive with constant $m=1$, it is easy to obtain the related results from \cite{MLY,WLX}, as well as the results from \cite[Corollary~3.1]{ZHJ} and \cite[Corollary~3.2]{ZHJ} \emph{without} assuming the convexity of the set $S$ therein. Note that it is not possible to apply \cite[Theorem~3.1]{ZHJ} to derive these results since the function $f(x):=J(x)+\delta(x;\O)$ never satisfies a center-Lipschitz condition at $\ox\in S_0$ if $\O$ is a proper subset of $X$.
\begin{Corollary}
Consider the perturbed distance function defined by {\rm (\ref{pd})} and the set $S_0$ given by {\rm (\ref{pds})} with $\ox\in S_0$. Suppose that $J$ satisfies a center-Lipschitz condition on $S$ with constant $\ell<1$. Then
\begin{equation*}
\Hat\partial d_\O^J(\ox)=\Hat\partial(J+\delta_\O)(\ox)\cap \B^*.
\end{equation*}
In particular, one has
\begin{equation*}
\Hat\partial d(\ox; \O)=\Hat N(\ox; \O)\cap \B^*,
\end{equation*}
where $\Hat N(\ox;\O):=\Hat\partial\delta_\O(\ox)$.
\end{Corollary}

\section{H\"{o}lder subdifferential formulas}

In this section, we develop H\"{o}lder subdifferential formulas for the infimal convolution (\ref{ic}). Given an extended-real-valued function $g: X\to \oR$ with $\ox\in \mbox{\rm dom}\,g$ and given $s>0$, the $s-$\emph{H\"{o}lder subdifferential} of $f$ at $\ox$ is defined by
\begin{equation*}
\partial_s g(\ox):=\big\{x^*\in X^*\; \big |\; \liminf_{x\to \ox}\dfrac{g(x)-g(\ox)-\la x^*, x-\ox\ra}{\|x-\ox\|^{1+s}}> -\infty\big\},
\end{equation*}
which reduces to the proximal subdifferential if $s=1$ and reduces to the subdifferential in the sense of convex analysis if $g$ is convex.

It follows from the definition that $x^*\in \partial_sg(\ox)$ if and only if there exist $\sigma>0$ and $\delta>0$ such that
\begin{equation*}
\la x^*, x-\ox\ra \leq g(x)-g(\ox)+\sigma \|x-\ox\|^{1+s}\; \mbox{\rm whenever }\|x-\ox\|<\delta.
\end{equation*}
\begin{Theorem}\label{t2} Consider the infimal convolution {\rm (\ref{ic})} in which $\ph(0)=0$ and consider the set $S_0$ given by {\rm (\ref{s})} with $\ox\in S_0$. The following hold:\\[1ex]
{\rm\bf (i)} $\partial_s T^f_\ph(\ox)\subset \partial_sf(\ox)\cap [-\partial_s \ph(0)].$\\
{\rm\bf (ii)} Suppose that $\ph$ is coercive on $X$ with constant $m>0$, and that $f$ satisfies a center-Lipschitz condition on $D:=\mbox{\rm dom}\,f$ with constant $\ell$ where $0\leq \ell <m$. Then
\begin{equation}\label{f2}
\partial_s T_\ph^f(\ox)= \partial_s f(\ox)\cap \big[-\partial_s \ph(0)\big].
\end{equation}
\end{Theorem}
{\bf Proof.} {\bf (i)} Fix any $x^*\in \partial_s T^f_\ph(\ox)$. Then there exist $\sigma>0$ and $\delta>0$ such that
\begin{equation*}
\la x^*, x-\ox\ra \leq T_\ph^f(x)-T_\ph^f(\ox)+\sigma \|x-\ox\|^{1+s}\; \mbox{\rm whenever }\|x-\ox\|<\delta.
\end{equation*}
Since $T_\ph^f(x)\leq f(x)$ for every $x\in X$, and $T_\ph^f(\ox)=f(\ox)$ as $\ox\in S_0$, one has
\begin{equation*}
\la x^*, x-\ox\ra \leq f(x)-f(\ox)+\sigma \|x-\ox\|^{1+s}\; \mbox{\rm whenever }\|x-\ox\|<\delta.
\end{equation*}
This implies $x^*\in \partial_sf(\ox)$.

Fixing any $v\in X$ with $\|v\|<\delta$ yields $\ox-v\in \B(\ox; \delta)$, and hence
\begin{align*}
\la x^*, -v\ra &\leq T_\ph^f(\ox-v)-T_\ph^f(\ox)+\sigma \|v\|^{1+s}\\
&\leq f(\ox)+\ph(v)-f(\ox)+\sigma \|v\|^{1+s}=\ph(v)+\sigma \|v\|^{1+s},
\end{align*}
which implies $x^*\in -\partial_s \ph(0)$. We have proved {\bf (i)}.\\[1ex]
{\bf (ii)} Using {\bf (i)}, it remains to prove the opposite inclusion in {\bf (ii)}. Assume by contradiction that there exists an element $x^*\in X^*$ that belong to $\partial_s f(\ox)\cap \big[-\partial_s \ph(0)\big]$, but $x^*\notin \partial_sT^f_\ph(\ox)$. Then there exist sequences $(x_k)$ and $(\sigma_k)$ with $x_k\to \ox$ and $\sigma_k\to \infty$ as $k\to \infty$ and
\begin{equation}\label{con}
\la x^*, x_k-\ox\ra > T_\ph^f(x_k)-T_\ph^f(\ox)+\sigma_k \|x_k-\ox\|^{1+s}\; \mbox{\rm for every } k\in \N.
\end{equation}
This implies $x_k\neq \ox$ and
\begin{equation*}
T_\ph^f(x_k)\leq \|x^*\| \|x_k-\ox\|+f(\ox)\; \mbox{\rm for every }k\in \N.
\end{equation*}
Without loss of generality, we can assume that $\|x_k-\ox\|<1$ for every $k$. Fix a sequence $(y_k)$ in $X$ such that
\begin{equation*}
\ph(y_k-x_k)+f(y_k)< T_\ph^f(x_k)+\|x_k-\ox\|^{1+s}\leq \|x^*\| \|x_k-\ox\|+f(\ox)+\|x_k-\ox\|.
\end{equation*}
This implies $y_k\in \mbox{\rm dom}\,f$, and hence
\begin{align*}
m\|y_k-\ox\|&\leq m\|y_k-x_k\|+m\|x_k-\ox\|\\
&\leq \ph(y_k-x_k)+m\|x_k-\ox\|\\
&\leq f(\ox)-f(y_k)+\|x^*\|\|x_k-\ox\|+\|x_k-\ox\|+m\|x_k-\ox\|\\
&\leq \ell \|\ox-y_k\|+(m+\|x^*\|+1)\|x_k-\ox\|.
\end{align*}
It follows that
\begin{equation*}
\|y_k-\ox\|\leq (m-\ell)^{-1}(m+\|x^*\|+1)\|x_k-\ox\|=c\|x_k-\ox\|\to 0\; \mbox{\rm as }k\to \infty,
\end{equation*}
where $c:=(m-\ell)^{-1}(m+\|x^*\|+1)$. Since $x^*\in \partial_s f(\ox)$, we can choose $\sigma>0$ and $\delta>0$ such that
\begin{equation*}
\la x^*, x-\ox\ra \leq f(x)-f(\ox)+\sigma \|x-\ox\|^{1+s}\; \mbox{\rm whenever }\|x-\ox\|<\delta.
\end{equation*}
Thus we can assume without loss of generality that
\begin{equation*}
\la x^*, y_k-\ox\ra \leq f(y_k)-f(\ox)+\sigma \|y_k-\ox\|^{1+s}\; \mbox{\rm for every }k\in \N.
\end{equation*}
Using also the fact that $x^*\in -\partial_s \ph(0)$, one finds $\delta_1>0$ such that
\begin{equation*}
\la -x^*, u\ra \leq \ph(u)+\sigma \|u\|^{1+s}\; \mbox{\rm whenever }\|u\|<\delta_1.
\end{equation*}
Since $\|y_k-x_k\|\to 0$, we can assume loss of generality that
\begin{equation*}
\la -x^*, y_k-x_k\ra\leq \ph(y_k-x_k)+\sigma\|y_k-x_k\|^{1+s}\; \mbox{\rm for every }k,
\end{equation*}
which implies $\la -x^*, y_k-x_k\ra\leq \ph(y_k-x_k)+\sigma(c+1)^{1+s}\|x_k-\ox\|^{1+s}$ for every $k$. It follows that
\begin{align*}
T_\ph^f(x_k)-T_\ph^f(\ox)-\la x^*, x_k-\ox\ra &=T_\ph^f(x_k)-f(\ox)-\la x^*, x_k-y_k\ra -\la x^*, y_k-\ox\ra\\
&\geq T_\ph^f(x_k)-f(\ox)-\la x^*, x_k-y_k\ra-\big[f(y_k)-f(\ox)+\sigma \|y_k-\ox\|^{1+s}\big]\\
&= T_\ph^f(x_k)-f(y_k)-\la x^*, x_k-y_k\ra-\sigma \|y_k-\ox\|^{1+s}\\
&\geq \ph (y_k-x_k)-\|x_k-\ox\|^{1+s}-\la x^*, x_k-y_k\ra-\sigma \|y_k-\ox\|^{1+s}\\
&\geq -\|x_k-\ox\|^{1+s}-\sigma(c+1)^{1+s}\|x_k-\ox\|^{1+s}-\sigma \|y_k-\ox\|^{1+s}\\
&\geq -\|x_k-\ox\|^{1+s} -\sigma(c+1)^{1+s}\|x_k-\ox\|^{1+s}-\sigma c^{1+s}\|x_k-\ox\|^{1+s}\\
&=-\big[1+\sigma(c+1)^{1+s}+\sigma c^{1+s}\big]\|x_k-\ox\|^{1+s}
\end{align*}
Denoting $\gamma:=\big[1+\sigma(c+1)^{1+s}+\sigma c^{1+s}\big]$ yields
\begin{align*}
T_\ph^f(x_k)-T_\ph^f(\ox)-\la x^*, x_k-\ox\ra \geq -\gamma\|x_k-\ox\|^{1+s}\; \mbox{\rm for every }k\in \N,
\end{align*}
which is a contradiction to (\ref{con}). The proof is now complete. $\h$

Let us obtain the following immediate corollaries.
\begin{Corollary}\label{c2}
Consider the infimal convolution {\rm (\ref{pm})} and the set $S_0$ given by {\rm (\ref{ss})} with $\ox\in S_0$. Suppose that $f$ satisfies a center-Lipschitz condition on $D:=\mbox{\rm dom}\,f$ with constant $\ell$ where $0\leq \ell <\|F\|^{-1}$. Then
\begin{equation}\label{f2}
\partial_s T^f_F(\ox)= \partial_s f(\ox)\cap \big\{x^*\in X^*\; \big |\; \sup_{u\in F}\la -x^*, u\ra\leq 1\big\}.
\end{equation}
\end{Corollary}

\begin{Corollary}
Consider the perturbed distance function defined by {\rm (\ref{pd})} and the set $S_0$ given by {\rm (\ref{pds})} with $\ox\in S_0$. Suppose that $J$ satisfies a center-Lipschitz condition on $S$ with constant $\ell<1$. Then
\begin{equation*}
\partial_sd_\O^J=\partial_s(J+\delta_\O)(\ox)\cap \B^*.
\end{equation*}
In particular, one has
\begin{equation*}
\partial_s d(\ox; \O)=N_s(\ox; \O)\cap \B^*,
\end{equation*}
where $N_s(\ox;\O):=\partial_s\delta_\O(\ox)$.
\end{Corollary}
\begin{Remark}{\rm Following the procedure developed in \cite{MLY}, it is possible to obtain results on Fr\'echet subdifferentials of the infimal convolution (\ref{ic}) when the reference point does not belong to the set $S_0$ as well as related results for \emph{limiting/Mordukhovich subdifferentials} of this class of functions}
\end{Remark}


\begin{thebibliography}{99}

\bibitem{BFQ} Burke, J.V., Ferris, M.C., Qian, M.: On the Clarke subdifferential of the distance function of a closed
set. J. Math. Anal. Appl. \bf (166), 199--213 (1992).
\bibitem{BT} Bounkhel, M., Thibault, L.: On various notions of regularity of sets in nonsmooth analysis. Nonlinear
Anal. {\bf 48}, 223–246 (2002)

\bibitem{C} Clarke, F.H., Optimization and Nonsmooth Analysis, Wiley, New York (1983).
\bibitem{CLSW} Clarke, F.H., Ledyaev, Y.S., Stern, R.J., Wolenski, P.R., Nonsmooth Analysis
and Control Theory, Springer, New York (1998).

\bibitem{cowo} Colombo, G., Wolenski, P.R., The subgradient formula for the minimal time
function in the case of constant dynamics in Hilbert space. J.
Global Optim. \textbf{28}, 269–-282 (2004).

\bibitem{JH} Jiang, Y., He, Y.: Subdifferentials of a minimum time function in normed spaces. J. Math. Anal. Appl.
{\bf 358}, 410--418 (2009).

\bibitem{HN} He, Y., Ng, K.F.: Subdifferentials of a minimum time function in Banach spaces. J. Math. Anal. Appl.
{\bf 321}, 896--910 (2006).



\bibitem{MLY} Meng, L., Li, C., Yao, J.C.: Limiting subdifferentials of perturbed distance functions in Banach spaces.
Nonlinear Anal. {\bf 75}, 1483--1495 (2012).

\bibitem{mor} Mordukhovich, B.S., Variational Analysis and
Generalized Differentiation, I: Basic Theory, II: Applications,
Grundlehren Series (Fundamental Principles of Mathematical
Sciences), Vols. 330 and 331, Springer, Berlin (2006).

\bibitem{MN10} Mordukhovich, B.S., Nam, N.M., Limiting subgradients of minimal time functions
in Banach spaces, J. Global Optim. {\bf 46}, 615--633, (2010).

\bibitem{MN} Mordukhovich, B.S., Nam, N.M., Subgradients of minimal time functions
under minimal assumptions. J. Convex Anal. {\bf 18}, 915--947 (2011).


\bibitem{NZ} Nam, N. M., Zalinescu, C., Variational analysis of directional minimal time functions and applications to location problems, to appear in Set-Valued and Var. Anal. {\bf 21}, 405--430 (2013).

\bibitem{WLX} Wang, J.H., Li, C., Xu, H.K.: Subdifferentials of perturbed distance function in Banach spaces.
J. Global Optim. {\bf 46}, 489--501 (2010).

\bibitem{ZHJ} Zhang, Y., He, Y., Jiang, Y., Subdifferentials of a perturbed minimal time function in normed spaces, Optim Lett, in press.

\end{thebibliography}
\end{document}